\documentclass [12pt,a4paper,reqno]{amsart}
\input amssymb.sty

\usepackage{color}
\usepackage{array}
\usepackage{tabularx}
\usepackage{enumerate}
\usepackage{amsmath}
\usepackage{amsfonts}
\usepackage{amssymb}
\usepackage{bbm}

\newtheorem{definition}{Definition}[section]

\newtheorem{lemma}[definition]{Lemma}
\newtheorem{theorem}[definition]{Theorem}
\newtheorem{proposition}[definition]{Proposition}
\newtheorem{corollary}[definition]{Corollary}
\newtheorem{example}[definition]{Example}
\newtheorem{remark}[definition]{Remark}

\begin{document}

\newcommand{\fin}{$\Box$\\}
\newcommand{\ds}{\displaystyle}
\newcommand{\saut}[1]{\hfill\\[#1]}
\newcommand{\vsp}{\vspace{.15cm}}
\newcommand{\difrac}{\displaystyle \frac}
\newcommand{\dist}{\textrm{dist}}
\newcommand{\diam}{\mathrm{diam}\ }
\newcommand{\levy}{\mathscr{B}}
\newcommand{\sheet}{\mathbbm{B}}
\newcommand{\sifbm}{\mathbf{B}}
\newcommand{\alphar}{\texttt{\large $\boldsymbol{\alpha}$}}


\title[Stationarity and self-similarity characterization of the sifBm]{Stationarity and self-similarity characterization of the Set-indexed Fractional Brownian Motion}

\author{Erick Herbin}
\address{Ecole Centrale Paris, Grande Voie des Vignes, 92295 Chatenay-Malabry, France} \email{erick.herbin@gmail.com}
\author{Ely Merzbach}
\address{Dept. of Mathematics, Bar Ilan University, 52900 Ramat-Gan, Israel}\email{merzbach@macs.biu.ac.il}

\subjclass[2000]{62\,G\,05, 60\,G\,15, 60\,G\,17, 60\,G\,18}


\begin{abstract}
The set-indexed fractional Brownian motion (sifBm) has been defined by Herbin-Merzbach (2006a) for indices that are subsets of a metric measure space.
In this paper, the sifBm is proved to statisfy a strenghtened definition of increment stationarity. This new definition for stationarity property allows to get a complete characterization of this process by its fractal properties: The sifBm is the only set-indexed Gaussian process which is self-similar and has stationary increments.

Using the fact that the sifBm is the only set-indexed process whose projection on any increasing path is a one-dimensional fractional Brownian motion, the limitation of its definition for a self-similarity parameter $0<H<1/2$ is studied, as illustrated by some examples. When the indexing collection is totally ordered, the sifBm can be defined for $0<H<1$.
\end{abstract}


\maketitle



\section{Introduction}
In \cite{ehem}, the set-indexed fractional Brownian motion (sifBm) was defined among processes indexed by a collection of subsets of a measure metric space. The study of its properties showed fractal behaviour such as a kind of increment stationarity and self-similarity.
In addition, it is proved that the projection of a sifBm on an increasing path is a one-parameter time changed fractional Brownian motion. 
Fine properties of multi-dimensional parameter fractional Brownian motions are studied by several authors (see \cite{istas, tudor, local} for example).

In this paper, we extend the increment stationarity property defined in \cite{ehem}. Instead of considering a stationarity property on $\Delta X_C$ (for $C\in\mathcal{C}_0$) that only involves marginal distributions of the increment process, we consider a property of stationarity of the distribution of the whole process $\Delta X=\left\{\Delta X_C; C\in\mathcal{C}_0\right\}$. We obtain a strengthened definition for increment stationarity which is preserved under projections on flows (increasing paths). More precisely, we show that if $X$ is a set-indexed process satisfying this new property of stationarity, then its projection on any flow is a one-dimensional increment stationary process. For that reason, this new definition can be considered as the most natural one. The set-indexed fractional Brownian motion is proved to satisfy this property.

The new stationarity definition allows us to get the main result of this paper: a complete characterization of the set-indexed fractional Brownian motion as the only set-indexed mean-zero Gaussian process which satisfies the two properties of increment stationarity and self-similarity. This property thus extends the well-known characterization of one-parameter fractional Brownian motion.

The second point of this paper is the use of flows to understand the limitation of the general sifBm's definition for a parameter $H\in (0,1/2]$, as opposed to one-parameter fractional Brownian motion which is defined for $0<H<1$. In the latter case, the behaviour of the process leads to critical values for $H$ (see \cite{cheridito,chnu} for instance). 
Here we observe that the set-indexed fractional Brownian motion can be defined for $0<H<1$ when the indexing collection $\mathcal{A}$ is totally ordered. On the contrary, we give new examples of indexing collections $\mathcal{A}$ on which the sifBm cannot be defined for $H>1/2$.
\vspace{10pt}

The paper is organized as follows: in section 2, indexing collection and the set-indexed fractional Brownian motion are defined, and projections of the sifBm on flows are studied. Section \ref{charact} is devoted to the study of the sifBm along flows to get a deeper understanding of its properties. Among the results, we get a characterization of the sifBm by its projection along flows, which constitutes a converse of a result in \cite{ehem}.
We prove that a set-indexed process is a set-indexed fractional Brownian motion if and only if its projections on every flows are one-dimensional fractional Brownian motions.
This gives a good justification of the definition of the sifBm
and opens the door to a variety of applications. In \cite{cras}, a part of this result was presented and also an integral representation for the sifBm was given.

In section \ref{sectHhalf}, we use the flows to get a better understanding of the limitation $H\in (0,1/2]$. This fact was already observed for the fBm indexed by the sphere of $\mathbf{R}^N$ (see \cite{istas}). The examples given explain why the sifBm cannot be defined in general for $H>1/2$. As a byproduct, we prove that the cardinality of a totally ordered indexing collection cannot exceed the continuum.

In section \ref{sectfrac}, the new strengthened definition for increment stationarity of a set-indexed process is studied.

Then in section \ref{sectcharact}, we establish the fractal characterization of the set-indexed fractional Brownian motion.

\section{Projections of the sifBm on flows}

We follow \cite{ehem} for the framework and notation. Our processes are indexed
by an {\em indexing collection} $\mathcal{A}$ of compact subsets of a locally compact
metric space $\mathcal{T}$ equipped with a Radon measure $m$ (denoted $(\mathcal{T},m)$).

Let $\mathcal A(u)$ denotes the class of finite unions from sets belonging to $\mathcal A$.

\begin{definition}[Indexing collection]\label{basic}
A nonempty class $\mathcal{A}$ of compact, connected subsets of $\mathcal{T}$ is 
called an {\em indexing collection}
if it satisfies the following:
\begin{enumerate}
 \item $\emptyset\in\mathcal{A}$, and $A^{\circ}\neq A$ 
if $A\notin\left\{ \emptyset, \mathcal{T} \right\}$.
In addition, there exists an increasing sequence 
$\left(B_n\right)_{n\in\mathbf{N}}$ of sets in $\mathcal{A}(u)$ such that
$\mathcal{T}=\bigcup_{n\in\mathbf{N}}B^{\circ}_n$.

\item $\mathcal{A}$ is closed under arbitrary intersections and if 
$A,B\in\mathcal{A}$ are nonempty, then $A\cap B$ is nonempty. 
If $(A_i)$ is an increasing sequence in $\mathcal{A}$ and if there exists
$n\in\mathbf{N}$ s. t. for all $i$, $A_i\subseteq B_n$ 
then $\overline{\bigcup_i A_i}\in\mathcal{A}$.

\item The $\sigma$-algebra generated by $\mathcal{A}$, $\sigma (\mathcal{A})=\mathcal{B}$,
the collection of all Borel sets of $\mathcal{T}$.

\item {\em Separability from above}\\
There exists an increasing sequence of finite subclasses 
$\mathcal{A}_n=\{A_1^n,...,A_{k_n}^n\}$ of $\mathcal{A}$ closed under intersections 
and satisfying $\emptyset, B_n \in\mathcal{A}_n(u)$  
and a sequence of functions 
$g_n:\mathcal{A}\rightarrow\mathcal{A}_n(u)\cup\left\{\mathcal{T}\right\}$ 
such that
\begin{enumerate}
\item $g_n$ preserves arbitrary intersections and finite unions \\
(i. e. $g_n(\bigcap_{A\in\mathcal{A}'}A )=\bigcap_{A\in\mathcal{A}'}g_n(A )$ 
for any $\mathcal{A}'\subseteq\mathcal{A}$, and \\
if $\bigcup_{i=1}^kA_i=\bigcup_{j=1}^mA_j'$, then
$\bigcup_{i=1}^kg_n(A_i)=\bigcup_{j=1}^mg_n(A_j')$);

\item for each $A\in\mathcal{A}$, $A\subseteq (g_n(A))^{\circ}$;

\item $g_n(A)\subseteq g_m(A)$ if $n\geq m$;

\item for each $A\in\mathcal{A}$, $A=\bigcap_n g_n(A)$;

\item if $A,A'\in\mathcal{A}$ then for every $n$, $g_n(A) \cap A' \in\mathcal{A}$, and if 
$A'\in\mathcal{A} _n$ then $g_n(A) \cap A' \in\mathcal{A} _n$;
\item   $g_n(\emptyset)=\emptyset~\forall n$.
\end{enumerate}

\item Every countable intersection of sets in $\mathcal{A}(u)$ may be expressed as
the closure of a countable union of sets in $\mathcal{A}$.

\end{enumerate}
(Note: ` $\subset$' indicates strict inclusion and `$\overline{(\cdot)}$' 
and`$(\cdot )^{\circ}$' denote respectively the closure and the interior of a set.)
\end{definition}

The {\em set-indexed fractional Brownian motion (sifBm) on $(\mathcal{T},\mathcal{A},m)$} was defined as the
centered Gaussian process
$\sifbm^H=\left\{\sifbm^H_U;\;U\in\mathcal{A}\right\}$ such that
\begin{equation}\label{defsifbm}
\forall \ U,V\in\mathcal{A};\  E\left[\sifbm^H_U \sifbm^H_V\right]=\frac{1}{2}
\left[m(U)^{2H}+m(V)^{2H}-m(U\bigtriangleup V)^{2H}\right],
\end{equation}
where $0<H\leq\frac{1}{2}$.

If $\mathcal{A}$ is provided with a structure of group on $\mathcal{T}$, properties of increment stationarity and self-similarity are studied in \cite{ehem}.
In the special case of $\mathcal{A}=\left\{ [0,t];\;t\in\mathbf{R}^N_+ \right\} \cup\left\{\emptyset\right\}$, we get a multiparameter process called Multiparameter fractional Brownian motion (MpfBm), whose properties are studied in \cite{mpfbm}.

The notion of flow is the key to reduce the proof of many theorems. It was
extensively studied in \cite{cime} and \cite{Ivanoff}.

\begin{definition}\label{flowdef}
An {\em elementary flow} is defined to be a continuous increasing function $f:[a,b]\subset\mathbf{R}_+\rightarrow\mathcal{A}$, i. e. such that
\begin{align*}
\forall s,t\in [a,b];\quad & s<t \Rightarrow f(s)\subseteq f(t)\\
\forall s\in [a,b);\quad & f(s)=\bigcap_{v>s}f(v)\\
\forall s\in (a,b);\quad & f(s)=\overline{\bigcup_{u<s}f(u)}.
\end{align*}

A {\em simple flow} is a continuous function $f:[a,b]\rightarrow\mathcal{A}(u)$ such that there exists a finite sequence $(t_0,t_1,\dots,t_n)$ with $a=t_0<t_1<\dots<t_n=b$ and elementary flows $f_i:[t_{i-1},t_i]\rightarrow\mathcal{A}$ ($i=1,\dots,n$) such that
\begin{equation*}
\forall s\in [t_{i-1},t_i];\quad
f(s)=f_i(s)\cup \bigcup_{j=1}^{i-1}f_j(t_j).
\end{equation*}
The set of all simple (resp. elementary) flows is denoted by $S(\mathcal{A})$ (resp. $S^e(\mathcal{A})$).
\end{definition}

In \cite{Ivanoff}, the projection of a set-indexed process $X=\left\{X_U;\;U\in\mathcal{A}\right\}$ on any elementary flow $f$ was considered as the real-parameter process 
$X^f=\left\{X_{f(t)};\;t\in [a,b]\subset\mathbf{R}_+\right\}$. Here, we define another parametrization of this projection, which allows simpler statements in the sequel.

\begin{definition}
For any set-indexed process $X=\left\{X_U;\;U\in\mathcal{A}\right\}$ on the space $(\mathcal{T},\mathcal{A},m)$ and any elementary flow $f:[a,b]\rightarrow\mathcal{A}$, we define the {\em $m$-standard projection of $X$ on $f$} as the process $$X^{f,m}=\left\{X^{f,m}_t=X_{f\circ\theta^{-1}(t)};\;t\in [a,b]\right\},$$
where $\theta:t\mapsto m[f(t)]$ and $\theta^{-1}$ is its pseudo-inverse function.
\end{definition}

The use of this new notation $X^{f,m}$ avoids any confusion with the projection $X^f$ previously defined.

Notice that since $\theta$ is non-decreasing, the function $\theta^{-1}$ is well-defined and  for all $t\in [a,b]$, we have $\theta(\theta^{-1}(t))=t$.

The following result, proved in \cite{ehem}, gives a good justification of the
definition of the sifBm.

\begin{proposition}\label{prop1}
Let $\sifbm^H$ be a sifBm on $(\mathcal{T},\mathcal{A},m)$ and $f:[a,b]\rightarrow\mathcal{A}$ be an elementary flow. Then the process
$(\sifbm^H)^{f,m}=\{\sifbm_{f\circ\theta^{-1}(t)}^H,\ t\in [a,b]\}$, where $\theta:t\mapsto m[f(t)]$, is a one-parameter fractional Brownian motion.
\end{proposition}

In section \ref{charact}, we prove the converse to Proposition \ref{prop1}. For
this purpose, we will use the following lemma proved in \cite{cime}.

\begin{lemma}\label{lemprojflow}
The finite dimensional distributions of an additive $\mathcal
A$-indexed process $X$ determine and are determined by the finite dimensional
distributions of the class $\{X^f,\ f\in S(\mathcal A)\}$.
\end{lemma}

\section{Characterization of the sifBm}\label{charact}

In the case of $L^2$-monotone outer-continuous set-indexed processes, we prove that the sifBm could be defined as a process whose projections on elementary flows is a one-dimensional fractional Brownian motion.

Recall the following definition (see \cite{Ivanoff})
\begin{definition}
A set-indexed process $X=\left\{X_U;\;U\in\mathcal{A}\right\}$ is said $L^2$-monotone outer-conti-nuous if for any decreasing sequence $\left(U_n\right)_{n\in\mathbf{N}}$ of sets in $\mathcal{A}$, 
\begin{equation*}
E\left[|X_{U_n}-X_{\bigcap_{k\in\mathbf{N}} U_k}|^2\right]\rightarrow 0
\end{equation*}
as $n\rightarrow\infty$.
\end{definition}

\begin{proposition}
The sifBm $\sifbm^H$ ($0<H\leq 1/2$) is $L^2$-monotone outer-continuous.
\end{proposition}

\proof
Let $\left(U_n\right)_{n\in\mathbf{N}}$ be a decreasing sequence in $\mathcal{A}$.
As $\bigcap_{k\in\mathbf{N}} U_k \in\mathcal{A}$, by definition of sifBm, we have
\begin{equation*}
\forall n\in\mathbf{N};\quad
E\left[|\sifbm^H_{U_n}-\sifbm^H_{\bigcap_{k\in\mathbf{N}} U_k}|^2\right]=
m(U_n\setminus \bigcap_{k\in\mathbf{N}} U_k)^{2H}
\end{equation*}
But, as $\left(U_n\right)_{n\in\mathbf{N}}$ is decreasing, by definition of a measure,
\begin{equation*}
m(U_n\setminus \bigcap_{k\in\mathbf{N}} U_k) \rightarrow 0.
\end{equation*}
Then the result follows.
\fin

The following lemma will be useful for the converse of proposition \ref{prop1}, and will be strenghtened in section \ref{sectHhalf} to understand links between structure of $\mathcal{A}$ and flows.

\begin{lemma}\label{lemflow}
For any $U_1,U_2,\dots,U_n\in\mathcal{A}$ such that $U_i\subset U_{i+1}$ ($\forall i=1,\dots,n-1$), there exist an elementary flow $f:\mathbf{R}_+\rightarrow\mathcal{A}$ and real numbers $0<t_1<t_2<\dots<t_n$ such that
\begin{equation*}
\forall i=1,\dots,n;\quad
f(t_i)=U_i.
\end{equation*}
\end{lemma}

\proof
This result is a particular case of  lemma 5.1.7 in \cite{Ivanoff} (and lemma 5 in \cite{cime}). As the sequence $U_1,U_2,\dots,U_n$ is increasing, $\mathcal{A}'=\left\{U_1,U_2,\dots,U_n\right\}$ constitutes a semilattice of $\mathcal{A}$ with a consistent numbering. The proof of lemma 5.1.7 in \cite{Ivanoff} constructs such an elementary flow $f$. Here the increasing property of $(U_i)_{1\leq i\leq n}$ allows $f$ to take its values in $\mathcal{A}$ ($\subset\mathcal{A}(u)$).
\fin

\begin{theorem}\label{thcharactflows}
Let $X=\left\{ X_U ;\;U\in\mathcal{A} \right\}$ be an $L^2$-monotone outer-continuous set-indexed process.

If the projection $X^f$ of $X$ on any elementary flow $f$, is a time-changed one-parameter fractional Brownian motion of parameter $H\in (0,1/2]$,
then there exists a Borel measure $\nu$ on $\mathcal{T}$ such that $X$ is a set-indexed fractional Brownian motion on $(\mathcal{T},\mathcal{A},\nu)$.
\end{theorem}

This theorem states that the time-changes giving to projections the law of a one-parameter fBm, determine a Borel measure $\nu$ such that $X$ is a sifBm on the space $(\mathcal{T},\mathcal{A},\nu)$.

A sketch of the proof is given in \cite{cras}. Here we present a complete proof. In particular, the importance of lemma \ref{lemflow} is pointed out.

\proof
Let $f:[a,b]\rightarrow\mathcal{A}$ be an elementary flow. As the projected process $X^f$ is a time-changed fBm of parameter $H$, we have
\begin{equation}\label{flowfbm}
\forall s,t\in [a,b];\quad
E\left[X^f_t - X^f_s\right]^2=|\theta_f(t) - \theta_f(s)|^{2H}
\end{equation}
where $\theta_f:\mathbf{R}_+\rightarrow\mathbf{R}_+$ is an increasing function.

The idea of the proof is the construction of a measure $\nu$ such that for any $f\in S^e(\mathcal{A})$,
\begin{equation*}
\forall t\in [a,b];\quad \theta_f(t)=\nu\left[f(t)\right].
\end{equation*}

For all $U\in\mathcal{A}$, let us define
\begin{equation*}
F^e_U = \left\{ f \in S^e(\mathcal{A}) :\; \exists u_f\in [a,b]; U= f(u_f) \right\}.
\end{equation*} 
As for all $f$ and $g$ in $F^e_U$, 
$\theta_f(u_f)^{2H} = \theta_g(u_g)^{2H} = E\left[X_U\right]^2$,
one can define
\begin{equation}\label{defpsi}
\psi(U)=\theta_f(u_f)=\left(E\left[X_U\right]^2\right)^{\frac{1}{2H}}.
\end{equation}
For all $U$ and $V$ in $\mathcal{A}$ with $U\subset V$, lemma \ref{lemflow} implies the existence of an elementary flow $f$ such that
\begin{equation*}
\exists u_f,v_f\in [a,b];\; u_f \leq v_f;\quad
U=f(u_f) \subset f(v_f)=V
\end{equation*}
Then, as the time-change $\theta_f$ is increasing, $\psi$ is non-decreasing in $\mathcal{A}$.

The definition of $\psi$ on $\mathcal{A}$ can be extended to $C=U\setminus\bigcup_{1\leq i \leq n}U_i$ by the inclusion-exclusion formula
\begin{align}\label{psiC}
\psi(C)=\psi(U)-\sum_{i=1}^n \psi\left(U\cap U_i\right)
&+\sum_{i<j}\psi\left(U\cap (U_i\cap U_j)\right) \nonumber\\
&-\dots+(-1)^n \psi\left(U\cap \left(\bigcap_{1\leq i\leq n}U_i\right)\right).
\end{align} 
We denote this class of sets by $\mathcal{C}$.

The definition (\ref{psiC}) of $\psi$ can be easily extended to the set $\mathcal{C}(u)$ of finite unions of elements of $\mathcal{C}$ in the same way.

A direct consequence of definition (\ref{psiC}) is that, for all $C_1, C_2\in\mathcal{C}$ such that $C=C_1\cup C_2 \in\mathcal{C}$, 
\begin{equation}\label{psiunionC}
\psi(C_1\cup C_2) = \psi(C_1) + \psi(C_2) - \psi(C_1\cap C_2)
\end{equation}
Let us remark that equality (\ref{psiunionC}) holds for any $C_1, C_2\in\mathcal{C}(u)$.

From the pre-measure $\psi$ defined on $\mathcal{C}$, the function
\begin{equation}\label{defm}
\nu : E\subset\mathcal{T} \mapsto \inf_{\genfrac{}{}{0pt}{}{C_i\in\mathcal{C}}{E\subset\cup C_i}}
\sum_{i=1}^{\infty}\psi(C_i)
\end{equation}
defines an outer measure on $\mathcal{T}$ (see \cite{rogers} pp. 9--26). 
Let us show that $\nu$ defines a Borel measure on the topological space $\mathcal{T}$.

Let $\mathcal{M}_{\nu}$ be the $\sigma$-field of $\nu$-measurable subsets of $\mathcal{T}$. It is known that $\nu$ is a measure on $\mathcal{M}_{\nu}$ (see \cite{rogers}, thm. 3). By definition, any $U\in\mathcal{A}$ is $\nu$-measurable if 
\begin{equation*}
\forall A\subset U, \forall B\subset\mathcal{T}\setminus U;\quad
\nu(A\cup B)=\nu(A)+\nu(B)
\end{equation*}
As the inequality $\nu(A\cup B)\leq \nu(A)+\nu(B)$ follows from definition of any outer-measure, it remains to show the converse inequality.

Consider any sequence $\left(C_i\right)_{i\in\mathbf{N}}$ in $\mathcal{C}$ such that $A\cup B\subset\bigcup_{i}C_i$. 
The sequence $\left(C_i\right)_{i\in\mathbf{N}}$ can be decomposed in the elements $C_i,\ i\in I$ such that $C_i\cap U=\emptyset$ and the $C_i,\ i\in J$ such that $C_i\subset U$ (if $C_i\cap U\ne\emptyset$ and $C_i\not\subset U$, cut $C_i=C'_i\cup C''_i$ where $C'_i\subset U$ and $C''_i\cap U=\emptyset$).

As
\begin{equation*}
A\cup B \subset\left[\bigcup_{i\in I}C_i\right] \cup 
\left[\bigcup_{i\in J}C_i\right]
\end{equation*}
we get the implications
\begin{equation*}
\forall i\in I;\;C_i\cap U=\emptyset\quad\Rightarrow\quad
A\subset\bigcup_{i\in J}C_i
\end{equation*}
and
\begin{equation*}
\forall i\in J;\;C_i\subset U\quad\Rightarrow\quad
B\subset\bigcup_{i\in I}C_i.
\end{equation*}
Then,
\begin{align*}
\sum_{i=1}^{\infty}\psi(C_i) = 
\underbrace{\sum_{i\in I}\psi(C_i)}_{\geq \nu(B)}+
\underbrace{\sum_{i\in J}\psi(C_i)}_{\geq \nu(A)}
\end{align*}
which leads to $\nu(A\cup B)\geq \nu(A)+\nu(B)$.

We have proved that $\mathcal{A}\subset\mathcal{M}_{\nu}$. By definition of $\mathcal{A}$, the smallest $\sigma$-field containing $\mathcal{A}$ is the Borel $\sigma$-field $\mathcal{B}$. Therefore, $\mathcal{B}\subset\mathcal{M}_{\nu}$ and $\nu$ is a measure on $\mathcal{B}$.

The second part of the proof is to show that the measure $\nu$ is an extension of $\psi$, i.~e. 
\begin{equation}\label{ext}
\forall U\in\mathcal{A};\quad \nu(U)=\psi(U).
\end{equation}

\begin{itemize}

\item For any $U\in\mathcal{A}$, by definition of $\nu(U)$,
\begin{equation}\label{ext1}
\nu(U)=\inf_{\genfrac{}{}{0pt}{}{C_i\in\mathcal{C}}{U\subset\cup C_i}}
\sum_{i=1}^{\infty}\psi(C_i)
\leq \psi(U)
\end{equation}

\item To prove the converse inequality, consider $U\in\mathcal{A}$ and a sequence $\left(C_i\right)_{i\in\mathbf{N}}$ in $\mathcal{C}$ such that $U\subset\bigcup_{i}C_i$. For all $n\in\mathbf{N}^*$, we have 
\begin{equation*}
U\subset\bigcup_{1\leq i\leq n}C_i \cup \left[U\setminus\bigcup_{1\leq i\leq n}C_i\right].
\end{equation*}
Then, (\ref{psiunionC}) implies
\begin{align}\label{majpsiC}
\psi(U)&=\psi\left(\bigcup_{1\leq i\leq n}C_i\right)
+\psi\left(U\setminus\bigcup_{1\leq i\leq n}C_i\right)\nonumber\\
&\leq\sum_{i=1}^{\infty}\psi(C_i) + 
\psi\left(U\setminus\bigcup_{1\leq i\leq n}C_i\right).
\end{align}
Using $L^2$-monotone outer continuity of $X$ and proposition 1.4.8 in \cite{Ivanoff}, we have
\begin{equation}\label{convpsiC}
\lim_{n\rightarrow\infty}\psi\left(U\setminus\bigcup_{1\leq i\leq n}C_i\right)
=0
\end{equation}
Thus, (\ref{majpsiC}) and (\ref{convpsiC}) imply that for all sequence $\left(C_i\right)_{i\in\mathbf{N}}$ in $\mathcal{C}$ such that $U\subset\bigcup_{i}C_i$,
\begin{equation*}
\psi(U)\leq\sum_{i=1}^{\infty}\psi(C_i)
\end{equation*}
and then, by definition of $\nu(U)$
\begin{equation}\label{ext2}
\psi(U)\leq \nu(U)
\end{equation}

\end{itemize}
Equality (\ref{ext}) follows from (\ref{ext1}) and (\ref{ext2}).

From (\ref{defpsi}) and (\ref{ext}), the Borel measure $\nu$ defined by (\ref{defm}) satisfies
\begin{equation*}
\forall U\in\mathcal{A};\quad
E\left[X_U\right]^2=\psi(U)^{2H}=\nu(U)^{2H}.
\end{equation*}

Using the Borel measure $\nu$, consider a set-indexed fractional Brownian motion $Y$ on $(\mathcal{T},\mathcal{A},\nu)$ (which exists as $0<H\leq 1/2$), defined by
\begin{equation*}
\forall U,V\in\mathcal{A};\quad
E\left[Y_U Y_V\right]=\frac{1}{2} \left[ \nu(U)^{2H} + \nu(V)^{2H} - \nu(U\bigtriangleup V)^{2H} \right].
\end{equation*}
According to proposition 6.4 in \cite{ehem}, projections of $Y$ on any elementary flow $f:[a,b]\rightarrow\mathcal{A}$ is a time-change one-parameter fractional Brownian motion, i. e. such that
\begin{align*}
\forall s,t\in [a,b];\quad
E\left[Y^f_t - Y^f_s\right]^2 &=|\nu\left[f(t)\right] - \nu\left[f(s)\right]|^{2H}\\
&=|\theta_f(t) - \theta_f(s)|^{2H}
\end{align*}
Then, the projections of the set-indexed processes $X$ and $Y$ on any elementary flow have the same distribution. By additivity, this fact holds also on any simple flow. Thus, lemma \ref{lemprojflow} implies $X$ and $Y$ have the same law.
\fin

Considering only $m$-standard projections on flows, theorem \ref{thcharactflows} gives the following characterization of the sifBm.

\begin{corollary}
Let $X=\left\{ X_U ;\;U\in\mathcal{A} \right\}$ be an $L^2$-monotone outer-continuous set-indexed process.
The following two assertions are equivalent:
\begin{enumerate}[(i)]
\item for any elementary flow $f:[a,b]\rightarrow\mathcal{A}$, the $m$-standard projection of $X$ on $f$ is a one-parameter fractional Brownian motion of index $H\in (0,1/2]$;
\item $X$ is a set-indexed fractional Brownian motion of index $H\in (0,1/2]$ on $(\mathcal{T},\mathcal{A},m)$.
\end{enumerate}
\end{corollary}

\section{Can sifBm be defined for $H>1/2$?}\label{sectHhalf}
In \cite{ehem}, the set-indexed fractional Brownian motion is defined for a parameter $H\in (0,1/2]$. As one-dimensional fractional Brownian motion is defined for $H\in (0,1)$, a natural question arises: Are there conditions on the indexing collection $\mathcal{A}$ such that sifBm on $(\mathcal{T},\mathcal{A},m)$ can be defined for $H>1/2$? 
Projections on flows allow to answer this question.

Let $\Phi^H : \mathcal{A}\times\mathcal{A}\rightarrow\mathbf{R}$ denote the function
\begin{equation*}
\Phi^H : (U,V)\mapsto\frac{1}{2}\left\{m(U)^{2H}+m(V)^{2H}-m(U\bigtriangleup V)^{2H}\right\}.
\end{equation*}
The question is: In which cases $\Phi^H$ can be seen as the covariance function of a set-indexed process?
In the following, we can see that this question is related to the two different cases either $\mathcal{A}$ is totally ordered or not.

Let us first describe the particular structure of a totally ordered indexing collection.

\begin{proposition}\label{orderprop}
If the indexing collection $\mathcal{A}$ is totally ordered by the inclusion, then there exists a surjective elementary flow $f:\mathbf{R}_+\rightarrow\mathcal{A}$, i. e. such that 
\begin{equation*}
\forall U\in\mathcal{A};\quad
U\in f\left(\mathbf{R}_+\right).
\end{equation*}

\end{proposition}

\proof
By definition of an indexing collection, $\mathcal{A}$ can be discretized by the increasing sequence of finite subclasses $(\mathcal{A}_n)_{n\in\mathbf{N}}$. As subclasses $\mathcal{A}_n$ are finite and totally ordered, lemma \ref{lemflow} implies for all $n$, existence of an elementary flow $f_n:\mathbf{R}_+\rightarrow\mathcal{A}$ such that
\begin{equation*}
\mathcal{A}_n\subseteq f_n\left(\mathbf{R}_+\right).
\end{equation*}
Note that, by construction of flows $(f_n)$ (see \cite{Ivanoff}, lemma 5.1.7), we have 
\begin{equation}\label{succfn}
\forall t\in f_n^{-1}\left(\mathcal{A}_n\right),
\forall m\geq n;\quad f_m(t)=f_n(t).
\end{equation}
Let us define $\mathcal{I}$, the set of $s\in\mathbf{R}_+$ such that $f_m(s)\in\mathcal{A}_m$ for some $m\in\mathbf{N}$.\\
From the sequence $(f_n)_{n\in\mathbf{N}}$, we define the function $f:\mathbf{R}_+\rightarrow\mathcal{A}$ in the following way:
\begin{itemize}
\item For all $t\in\mathcal{I}$, there exists $m\in\mathbf{N}$ such that $f_m(t)\in\mathcal{A}_m$. By (\ref{succfn}), the sequence $(f_n(t))_{n\geq m}$ is constant. We can set $f(t)=f_m(t)$.

\item In the construction of lemma 5.1.6 in \cite{Ivanoff}, the subset $\mathcal{I}$ is proved to be dense. Let us define for all $t\notin\mathcal{I}$,
\begin{equation*}
f(t)=\bigcap_{\genfrac{}{}{0pt}{}{s\in\mathcal{I}}{s>t}}f(s).
\end{equation*}
\end{itemize}
Let us show that $f$ satisfies the conclusions of the proposition \ref{orderprop}.

\begin{itemize}
\item $f$ is non-decreasing:
for all $s,t\in\mathbf{R}_+$ such that $s<t$, we have clearly
\begin{align*}
\bigcap_{\genfrac{}{}{0pt}{}{u\in\mathcal{I}}{u>s}}f(u) \subseteq
\bigcap_{\genfrac{}{}{0pt}{}{u\in\mathcal{I}}{u>t}}f(u)
\end{align*}
and then, $f(s)\subseteq f(t)$.

\item $f$ passes through every elements of $\bigcup_{n\in\mathbf{N}}\mathcal{A}_n$: 
for any $U\in\bigcup_{n\in\mathbf{N}}\mathcal{A}_n$, there exist $m\in\mathbf{N}$ and $t_U\in\mathcal{I}$ such that $U=f_m(t_U)$. Then, by definition of $f$ on $\mathcal{I}$, we have $f(t_U)=U$.

\item $f$ is continuous: according to definition \ref{flowdef}, we must verify that $f(t)=\bigcap_{s>t}f(s)$ and $f(t)=\overline{\bigcup_{s<t}f(s)}$.
Using density of $\mathcal{I}$, the right-continuity of $f$ comes directly from its definition
\begin{align*}
f(t)=\bigcap_{\genfrac{}{}{0pt}{}{s\in\mathcal{I}}{s>t}}f(s)
=\bigcap_{s>t}f(s).
\end{align*}
For the second equality, in the proof of lemma 5.1.6 in \cite{Ivanoff}, it is proved that
\begin{align*}
\overline{\bigcup_{\genfrac{}{}{0pt}{}{s\in\mathcal{I}}{s<t}}f(s)}
=\bigcap_{\genfrac{}{}{0pt}{}{s\in\mathcal{I}}{s>t}}f(s).
\end{align*}
Then the density of $\mathcal{I}$ allows to conclude $f(t)=\overline{\bigcup_{s<t}f(s)}$.

\item $f$ passes through all the elements of $\mathcal{A}$: for all $U\in\mathcal{A}\setminus\bigcup_n\mathcal{A}_n$, its approximations satisfy
\begin{align*}
\forall n\in\mathbf{N};\quad
g_n(U)\in\mathcal{A}_n(u) = \mathcal{A}_n
\end{align*}
because $\mathcal{A}_n$ is totally ordered for all $n$.
Therefore, for all $n\in\mathbf{N}$, there exists $t_{g_n(U)}$ in $\mathbf{R}_+$ such that $g_n(U)=f_n(t_{g_n(U)})$,
and we can write
\begin{align*}
U=\bigcap_{n\in\mathbf{N}} g_n(U)
=\bigcap_{n\in\mathbf{N}} f_n(t_{g_n(U)})
=\bigcap_{n\in\mathbf{N}} f(t_{g_n(U)}).
\end{align*}
The sequence $(t_{g_n(U)})_{n\in\mathbf{N}}$ is non-increasing in $\mathbf{R}_+$:
By definition,
\begin{align*}
g_{n+1}(U)&=f_{n+1}(t_{g_{n+1}(U)}) \\
g_n(U)&=f_n(t_{g_n(U)})=f_{n+1}(t_{g_n(U)})
\end{align*}
using the construction of $(f_n)_{n\in\mathbf{N}}$.
Then, as $(g_n(U))_{n\in\mathbf{N}}$ is non-increasing and $f_{n+1}$ is non-decreasing, we get $t_{g_{n+1}(U)} \leq t_{g_{n}(U)}$ from $g_n(U)\subseteq g_{n+1}(U)$.

Consequently, $(t_{g_n(U)})_{n\in\mathbf{N}}$ converges to some $t_U\in\mathbf{R}_+$. 
Then, using the continuity of $f$, we get
\begin{align*}
\bigcap_{n\in\mathbf{N}}f(t_{g_n(U})=f(t_U)
\end{align*}
which proves that $U\in f\left(\mathbf{R}_+\right)$.
\end{itemize} 
\fin

In \cite{Ivanoff}, proposition 1.3.5 shows that by definition, an indexing collection is not allowed to be too big. More precisely, the cardinality of $\mathcal{A}$ cannot exceed the cardinality of $\mathcal{P}(\mathbf{R})$, the set of subsets of $\mathbf{R}$.
In the particular case of a totally ordered indexing collection, this upper bound for size of $\mathcal{A}$ can be sharpened. From surjectivity of the flow $f$ in proposition \ref{orderprop}, we can state

\begin{corollary}
If the indexing collection $\mathcal{A}$ is totally ordered by the inclusion, then its cardinality cannot exceed cardinality of $\mathbf{R}$.
\end{corollary}

\vspace{30pt}

From this study of the particular case of a totally ordered collection $\mathcal{A}$, we can prove the existence the set-indexed fractional Brownian motion for a parameter $H\in (0,1)$, as in one-parameter case.

\begin{theorem}\label{thH01}
If the indexing collection $\mathcal{A}$ is totally ordered by the inclusion, then the set-indexed fractional Brownian motion on $(\mathcal{T},\mathcal{A}, m)$ can be defined for $0<H<1$.
\end{theorem}

\proof
According to proposition \ref{orderprop}, as $\mathcal{A}$ is totally ordered, there exists an elementary flow $f:\mathbf{R}_+\rightarrow\mathcal{A}$ passing through every $U\in\mathcal{A}$, i. e.  such that 
\begin{align*}
\forall U\in\mathcal{A};\quad
U\in f\left(\mathbf{R}_+\right).
\end{align*}
Then, the existence of the sifBm is equivalent to the existence of its projection on the flow $f$. For any $H\in (0,1)$, let us consider a one-parameter fractional Brownian motion $B^H=\{B^H_t;\;t\in\mathbf{R}_+\}$ of self-similarity parameter $H$. \\
The set-indexed process 
$X=\left\{X_U=B^H_{\theta\circ f^{-1}(U)};\;U\in\mathcal{A}\right\}$,
where $\theta:t\mapsto m[f(t)]$, is a mean-zero Gaussian process such that for all $U, V\in\mathcal{A}$,
\begin{align*}
E\left[ X_UX_V \right]
&=\frac{1}{2}\left[ (\theta\circ f^{-1}(U))^{2H} + (\theta\circ f^{-1}(V))^{2H}
- |\theta\circ f^{-1}(U)-\theta\circ f^{-1}(V)|^{2H} \right] \\
&=\frac{1}{2}\left[ m(U)^{2H} + m(V)^{2H} - |m(U) - m(V)|^{2H} \right].
\end{align*}
As $\mathcal{A}$ is totally ordered, we have either $U\subseteq V$ or $V\subseteq U$, and then
\begin{equation*}
\left| m(U) - m(V) \right|=m(U\bigtriangleup V).
\end{equation*}
Thus, the covariance structure of $X$ is given by
\begin{equation*}
\forall U, V\in\mathcal{A};\quad
E\left[ X_U X_V \right]=\frac{1}{2}\left[ m(U)^{2H} + m(V)^{2H} - m(U\bigtriangleup V)^{2H} \right]
\end{equation*}
and it follows that $X$ is a sifBm of parameter $H$ on $(\mathcal{T},\mathcal{A},m)$.
\fin

Theorem \ref{thH01} suggests that if there exists a real limitation of sifBm's definition to $H<1/2$, it must be only for partially ordered indexing collections.
The following example shows that even in the simple case of rectangles in $\mathbf{R}^2$, the sifBm may not be defined for $H>1/2$.

\begin{example}
Let us consider the indexing collection constituted by rectangles of $\mathbf{R}^2_+$,
$\mathcal{A}=\left\{ [0,t]; t\in\mathbf{R}^2_+ \right\} \cup \left\{\emptyset\right\}$.
$\mathcal{A}$ is not totally ordered, and then, theorem \ref{thH01} cannot be applied.
In fact, for $H>1/2$ the function $\Phi^H$ is not positive definite.\\
Let us consider the four points of $\mathbf{R}^2_+$, defined by their coordinates $t_1=(1,1)$, $t_2=(2,1)$, $t_3=(1,2)$ and $t_4=(2,2)$. The points $t_i$ ($1\leq i\leq 4$) define four elements $U_i=[0,t_i]$ of $\mathcal{A}$. We compute $\Phi^H(U_i,U_j)$ for all $1\leq i,j \leq 4$ ($m$ is the Lebesgue measure in $\mathbf{R}^2$). The diagonal terms are
\begin{align*}
&\Phi^H(U_1,U_1) = m(U_1)^{2H} = 1; \quad
&\Phi^H(U_2,U_2) = m(U_2)^{2H} = 2^{2H};\\
&\Phi^H(U_3,U_3) = m(U_3)^{2H} = 2^{2H}; \quad
&\Phi^H(U_4,U_4) = m(U_4)^{2H} = 4^{2H}.
\end{align*}
The cross terms are
\begin{align*}
\Phi^H(U_1,U_2) &= \frac{1}{2}\left[m(U_1)^{2H}+m(U_2)^{2H}-m(U_2\setminus U_1)^{2H}\right] \\
&= \frac{1}{2}\left[1+2^{2H}-1\right]= 2^{2H-1};\\
\Phi^H(U_1,U_3) &= 2^{2H-1}; \\
\Phi^H(U_1,U_4) &= \frac{1}{2}\left[1+4^{2H}-3^{2H}\right];
\end{align*}
and
\begin{align*}
\Phi^H(U_2,U_3) &= \frac{1}{2}\left[2^{2H}+2^{2H}-2^{2H}\right]=2^{2H-1}; \\
\Phi^H(U_2,U_4) &= \frac{1}{2}\left[2^{2H}+4^{2H}-2^{2H}\right]=2^{4H-1}; \\
\Phi^H(U_3,U_4) &= 2^{4H-1}.
\end{align*}
By computation, the matrix 
\begin{equation*}
\left( \begin{array}{c c c c}
1 & 2^{2H-1} & 2^{2H-1} & \frac{1+2^{4H}-3^{2H}}{2} \\
2^{2H-1} & 2^{2H} & 2^{2H-1} & 2^{4H-1} \\
2^{2H-1} & 2^{2H-1} & 2^{2H} & 2^{4H-1} \\
\frac{1+2^{4H}-3^{2H}}{2} & 2^{4H-1} & 2^{4H-1} & 2^{4H}
\end{array} \right)
\end{equation*}
is not positive definite for $H=3/4$ (although it is for $H=1/2$).
Therefore $\Phi^{3/4}$ is not positive definite and consequently, the sifBm cannot be defined on $(\mathbf{R}^2_+, \mathcal{A}, m)$ for $H=3/4$.
\end{example}

The following example shows that sifBm's definition can be used to obtain an extension of fractional Brownian motion indexed by a differential manifold. In that case, the choice of the indexing collection on the manifold is crucial and can lead to different processes, whose definitions are limited to $0<H\leq 1/2$ or not.

\begin{example}\label{excercle}
Suppose we aim to extend fractional Brownian motion for indices in the unit circle $\mathbb{S}_1$ in $\mathbf{R}^2$. Let us fix a point $O\in\mathbb{S}_1$ and define $\mathcal{A}$ as the collection $\left\{\overset{\curvearrowright}{0M}; M\in\mathbb{S}_1\right\}\cup\left\{\emptyset\right\}$, where $\overset{\curvearrowright}{0M}$ denotes the positive oriented arc from $O$ to $M$. $\mathcal{A}$ is clearly an indexing collection which is totally ordered. Then, theorem \ref{thH01} implies the existence of a sifBm on $(\mathbb{S}_1, \mathcal{A}, m)$ for a parameter $H\in (0,1)$, where $m$ denotes the Lebesgue measure on $\mathbb{S}_1$. It is defined as the mean-zero Gaussian process $X=\left\{X_M ; M\in\mathbb{S}_1\right\}$ such that
\begin{equation*}
\forall M,M'\in\mathbb{S}_1;\quad
E\left[X_M-X_{M'}\right]^2 = m(\overset{\curvearrowright}{0M}\bigtriangleup\overset{\curvearrowright}{0M'})^{2H}=m(\overset{\curvearrowright}{MM'})^{2H}.
\end{equation*}

Another choice of indexing collection is $\mathcal{A}'=\left\{\overset{\frown}{0M}; M\in\mathbb{S}_1\right\}\cup\left\{\emptyset\right\}$, where $\overset{\frown}{0M}$ denotes the smallest arc of circle from $O$ to $M$. As $\mathcal{A}'$ is not totally ordered, there is a priori a limitation of sifBm's definition on $(\mathbb{S}_1, \mathcal{A}', m)$ for a parameter $H\in (0,1/2)$.

Another point of view is followed in Istas' extension of fractional Brownian motion indexed by points on the unit circle, considered as a metric space (and not as a measure space). In \cite{istas}, the {\em periodical fractional Brownian motion (PFBM)} is defined as the mean-zero Gaussian process $X=\left\{X_M ; M\in\mathbb{S}_1\right\}$ such that $X_O = 0$ (for some $O\in\mathbb{S}_1$) almost surely and
\begin{align*}
\forall M,M'\in\mathbb{S}_1;\quad
E\left[ X_M - X_{M'} \right]^2 = \left[d(M,M')\right]^{2H}
\end{align*}
where $d(M,M')$ is the distance between $M$ and $M'$ on the circle $\mathbb{S}_1$.

This process is different from the two previous definitions based on set-indexed fractional Brownian motion, in the sense that the covariance function cannot be expressed in terms of some measure on $\mathbb{S}_1$. 

However, if we only consider the positive half-circle $\frac{1}{2}\mathbb{S}_1$ starting from $O\in\mathbb{S}_1$, then
\begin{align*}
\forall M,M'\in\frac{1}{2}\mathbb{S}_1;\quad
m(\overset{\frown}{0M}\bigtriangleup\overset{\frown}{0M'})=
m(\overset{\curvearrowright}{MM'})^{2H} = \left[d(M,M')\right]^{2H}
\end{align*}
and the three covariance functions are identical.
Therefore the three different processes defined on $\mathbb{S}_1$ coincide on $\frac{1}{2}\mathbb{S}_1$. In that sense, Istas' PFBM on the half-circle can be seen as a particular case of the sifBm. Consequently, fractal properties such as stationarity and self-similarity are satisfied by this process on $\frac{1}{2}\mathbb{S}_1$ (cf. section \ref{sectfrac}) but stationarity cannot hold on the whole unit circle.

Moreover, as seen later, the characterization by fractal properties leads naturally to our first definition (cf. section \ref{sectcharact}).\end{example}

\section{Fractal properties}\label{sectfrac}
\subsection{Increment stationarity}

The increments of a set-indexed process are defined from the collection of subsets $\mathcal{C}$.\\
For all $C=U\setminus\bigcup_{1\leq i\leq n} U_i$, we define the increment of
the process $X$ on $C$ by
\begin{equation}
\Delta X_C=X_U-\sum_{i=1}^n X_{U\bigcap U_i}
+\sum_{i<j}X_{U\bigcap (U_i\bigcap U_j)} -\dots
+(-1)^n X_{U\bigcap \left(\bigcap_{1\leq i\leq n}U_i\right)}.
\end{equation}
This increment is always well defined for the sifBm $\sifbm^H$ since for all $U,V\in\mathcal{A}$ such that $U\cup V\in\mathcal{A}$, we have $E\left[X_U+X_V-X_{U\cap V}-X_{U\cup V}\right]^2=0$.

In \cite{ehem}, we defined a stationarity property for a set-indexed process $X=\left\{X_U;\;U\in\mathcal{A}\right\}$ on $(\mathcal{T},\mathcal{A},m)$ by
\begin{equation}\label{oldstat}
\forall C,C'\in\mathcal{C}_0;\quad
m(C)=m(C') \Rightarrow \Delta X_C \stackrel{(d)}{=} \Delta X_{C'}
\end{equation}

where $\mathcal{C}_0$ denotes the set of elements $U\setminus V$ with $U,V\in\mathcal{A}$.
Property (\ref{oldstat}) is called $\mathcal{C}_0$-stationarity.

As $\mathcal{C}_0$-stationarity only concerns marginal distributions of the increment process $\Delta X$, but not distribution of the process, the property is weaker than the classical increment stationarity property for one-parameter processes.

In that view, it seems judicious to strengthen the increment stationarity definition in such a way that projections of a increment stationary set-indexed process on any flow give increment stationary one-parameter processes.

\begin{definition}\label{defstat}
A set-indexed process $X=\left\{X_U;\;U\in\mathcal{A}\right\}$ is said to have $m$-stationary $\mathcal{C}_0$-increments if for any integer $n$, for all $V\in\mathcal{A}$ and for all increasing sequences $(U_i)_{1\leq i\leq n}$ and $(A_i)_{1\leq i\leq n}$ in $\mathcal{A}$,
\begin{equation*}
\forall i,\;m(U_i\setminus V)=m(A_i)
\quad\Rightarrow\quad
\left(\Delta X_{U_1\setminus V},\dots,\Delta X_{U_n\setminus V}\right) \stackrel{(d)}{=}
\left(\Delta X_{A_1},\dots,\Delta X_{A_n}\right).
\end{equation*}
\end{definition}

\vspace{10pt}

\begin{proposition}\label{propsifbmstat}
The sifBm has $m$-stationary $\mathcal{C}_0$-increments.
\end{proposition}

\proof
Let $X=\left\{X_U;\;U\in\mathcal{A}\right\}$ be a sifBm.
For any integer $n$, let us consider $V\in\mathcal{A}$ and increasing sequences $(U_i)_{1\leq i\leq n}$ and $(A_i)_{1\leq i\leq n}$ in $\mathcal{A}$ such that $m(U_i\setminus V)=m(A_i)$ ($\forall 1\leq i\leq n$). Without loss of generality, we can assume $V\subset U_i$ ($\forall i$). Let us compute for all $\lambda_1,\dots,\lambda_n$ in $\mathbf{R}$,
\begin{align*}
E\left[\lambda_1 \Delta X_{U_1\setminus V}+\dots+\lambda_n\Delta X_{U_n\setminus V}\right]^2
&=\sum_{i,j}\lambda_i\lambda_j E\left[\Delta X_{U_i\setminus V} \Delta X_{U_j\setminus V}\right]\\
&=\sum_{i,j}\lambda_i\lambda_j E\left[(X_{U_i}-X_V)(X_{U_j}-X_V)\right]\\
&=\sum_{i,j}\lambda_i\lambda_j \left(E\left[X_{U_i}X_{U_j}\right]-E\left[X_{U_i}X_V\right]
-E\left[X_V X_{U_j}\right]+E\left[X_V\right]^2\right)\\
&=\sum_{i,j}\lambda_i\lambda_j \left(m(U_i\bigtriangleup V)^{2H}+m(U_i\bigtriangleup V)^{2H}-m(U_i\bigtriangleup U_j)^{2H}\right).
\end{align*}
Assumptions on $(U_i)_i$ and $(A_i)_i$ imply $m(U_i\bigtriangleup V)=m(U_i\setminus V)=m(A_i)$, and as $(U_i)_i$ is increasing, $m(U_i\bigtriangleup U_j)=|m(U_i\setminus V)-m(U_j\setminus V)|=|m(A_i)-m(A_j)|$.
Then,  for all $\lambda_1,\dots,\lambda_n$ in $\mathbf{R}$
\begin{equation*}
E\left[\lambda_1 \Delta X_{U_1\setminus V}+\dots+\lambda_n\Delta X_{U_n\setminus V}\right]^2=E\left[\lambda_1 \Delta X_{A_1}+\dots+\lambda_n \Delta X_{A_n}\right]^2
\end{equation*}
and, as the process $\Delta X$ is centered Gaussian, the result follows.
\fin

\begin{example}
Following the notation of example \ref{excercle}, the sifBm defined on $(\mathbb{S}_1,\mathcal{A},m)$, which provides an extension of fractional Brownian motion indexed by points of the unit circle of $\mathbf{R}^2$, has $m$-stationary $\mathcal{C}_0$-increments. By definition, $\mathcal{C}_0$ consists of all elements $\overset{\curvearrowright}{MM'}$ where $M,M'\in\mathbb{S}_1$. Then this stationarity property states that the law of the process $\Delta X$ is invariant by translations along $\mathbb{S}_1$.
\end{example}

The following proposition shows that definition \ref{defstat} provides a natural extension of stationarity property for one-parameter processes. Then, it justifies this new definition for stationarity of set-indexed processes.

\begin{proposition}\label{statproj}
A set-indexed process $X=\left\{X_U;\;U\in\mathcal{A}\right\}$ has $m$-stationary $\mathcal{C}_0$-increments if and only if the $m$-standard projection of $X$ on any elementary flow $f:\mathbf{R}_+\rightarrow\mathcal{A}$ has stationary increments, i. e.
\begin{equation*}
\left\{X^{f,m}_{t+h}-X^{f,m}_{h};\;t\in\mathbf{R}_+\right\} \stackrel{(d)}{=}
\left\{X^{f,m}_{t}-X^{f,m}_{0};\;t\in\mathbf{R}_+\right\}
\end{equation*}
where $X^{f,m}=\left\{X_{f\circ\theta^{-1}(t)};\;t\in\mathbf{R}_+\right\}$ and $\theta:t\mapsto m[f(t)]$.
\end{proposition}

\proof
We prove the first implication. Assume that $X$ has $m$-stationary $\mathcal{C}_0$-increments and that $f$ is an elementary flow.\\
For all $t_1<t_2<\dots<t_n$ and $h$ in $\mathbf{R}_+$, consider for $1\leq i\leq n$, $U_i=f\circ\theta^{-1}(t_i+h)$, $V=f\circ\theta^{-1}(h)$ and
\begin{equation*}
C_i=U_i\setminus V=f\circ\theta^{-1}(t_i+h)\setminus f\circ\theta^{-1}(h)
\quad\textrm{and}\quad
A_i=f\circ\theta^{-1}(t_i).
\end{equation*}
As $(U_i)_{1\leq i\leq n}$ and $(A_i)_{1\leq i\leq n}$ are increasing and $V\subset U_i$ ($\forall i$), we have
\begin{align*}
\Delta X_{U_i\setminus V}&=X_{U_i} - X_V\\
&=X_{f\circ\theta^{-1}(t_i+h)} - X_{f\circ\theta^{-1}(h)}\\
&=X^{f,m}_{t_i+h} - X^{f,m}_{h}
\end{align*}
and
\begin{align*}
m(U_i\setminus V)&=m[f\circ\theta^{-1}(t_i+h)] - m[f\circ\theta^{-1}(h)]\\
&=\theta\circ\theta^{-1}(t_i+h) - \theta\circ\theta^{-1}(h)\\
&=t_i\\
&=m(A_i).
\end{align*}
Then, $m$-stationarity of the set-indexed process $X$ implies
\begin{equation*}
(\Delta X_{U_1\setminus V},\dots,\Delta X_{U_n\setminus V}) \stackrel{(d)}{=}
(X_{A_1},\dots,X_{A_n})
\end{equation*}
which gives
\begin{equation*}
(X^{f,m}_{t_1+h}-X^{f,m}_h,\dots,X^{f,m}_{t_n+h}-X^{f,m}_h) \stackrel{(d)}{=}
(X^{f,m}_{t_1},\dots,X^{f,m}_{t_n})
\end{equation*}
and the increment stationarity of the $m$-standard projection of $X$ on $f$.

In the same way, using lemma \ref{lemflow} to define a flow passing through every $U_i$ ($1\leq i\leq n$), we prove the converse.
\fin

\subsection{Self-similarity}

In \cite{ehem}, we defined the self-similarity property of a set-indexed process with respect to action of a group $G$ satisfying the following assumptions.

We suppose that $\mathcal{A}$ is provided with the operation of a non trivial 
group $G$ that can be extended satisfying 
\begin{align}\label{opG}
\forall U,V\in\mathcal{A}, \forall g\in G;\quad
&g.(U\cup V)=g.U \cup g.V\\
&g.(U\setminus V)=g.U \setminus g.V \nonumber
\end{align}
and assume
there exists a surjective function $\mu : G\rightarrow\mathbf{R}^{*}_{+}$
\begin{equation}\label{mu}
\forall U\in\mathcal{A}, \forall g\in G;\quad
m(g.U)=\mu(g).m(U) .
\end{equation}


A set-indexed process $X=\left\{X_U;\; U\in\mathcal{A}\right\}$ is said to be
self-similar of index $H$, if there exists a group $G$ which operates on
$\mathcal{A}$, and satisfies (\ref{opG}) and (\ref{mu}),
such that for all $g\in G$,
\begin{equation}
\left\{X_{g.U};\; U\in\mathcal{A} \right\}
\stackrel{(d)}{=} \left\{\mu(g)^H .X_U;\; U\in\mathcal{A} \right\}
\end{equation}

\begin{remark}\label{remRN}
In the case of $\mathcal{A}=\left\{ [0,u];\;u\in\mathbf{R}^N_+\right\}\cup\{\emptyset\}$, the operation of $G=\mathbf{R}_+$ defined by
\begin{align*}
\mathbf{R}_+ \times \mathcal{A} &\rightarrow \mathcal{A}\\
(a,[0,u]) &\mapsto [0,au]
\end{align*}
satisfies assumptions (\ref{opG}) and (\ref{mu}). 
\end{remark}

On the contrary to stationarity property, self-similarity of a set-indexed process does not imply self-similarity of projections on flows in a natural way.
This is essentially due to the fact that there is no connection between zooming in $\mathcal{A}$ (through operation of $G$) and zooming along a flow (through multiplication by $\mathbf{R}_+$).
For instance,  in the frame of remark \ref{remRN}, if $[0,u]\in\mathcal{A}$ belongs to an elementary flow $f$, i. e. if there exists $t\in\mathbf{R}_+$ such that $f(t)=[0,u]$, the element $[0,au]$ ($a>0$) of $\mathcal{A}$ does not belong necessarily to $f$.

However, under some additional assumptions either about flows or about the set-indexed process, standard projections can inherit the self-similarity property. 

\begin{proposition}\label{propself}
Let $X=\left\{X_U;\; U\in\mathcal{A}\right\}$ be a set-indexed process on $(\mathcal{T},\mathcal{A},m)$ which satisfies the two following properties:
\begin{enumerate}
\item self-similarity of index $H$ (with respect to operation of a group $G$ satisfying assumptions (\ref{opG}) and (\ref{mu})), 
\item $m$-stationarity of  $\mathcal{C}_0$-increments.
\end{enumerate}
Then, the $m$-standard projection of $X$ on any elementary flow $f$ is self-similar of index $H$, i.~ e.
\begin{equation*}
\forall a\in\mathbf{R}_+;\quad
\left\{X^{f,m}_{at};\; t\in\mathbf{R}_+ \right\}
\stackrel{(d)}{=} \left\{a^H .X^{f,m}_t;\; t\in\mathbf{R}_+ \right\}
\end{equation*}
where $X^{f,m}_t=X_{f\circ\theta^{-1}(t)}$ and $\theta:t\mapsto m[f(t)]$.
\end{proposition}

\proof
Let $f$ be any elementary flow, $a\in\mathbf{R}_+$ and $t_1<t_2<\dots<t_n$ a sequence of elements of $\mathbf{R}_+$.
For all $i=1,\dots,n$, consider $U_i=f\circ\theta^{-1}(t_i)$.\\
As $\mu$ is a surjective function, there exists $g\in G$ such that $a=\mu(g)$.\\
As
\begin{align*}
\forall i=1,\dots,n;\quad
m(f\circ\theta^{-1}(at_i))&=\theta\circ\theta^{-1}(at_i)=at_i\\
&=\mu(g) m(U_i)=m(g.U_i),
\end{align*}
by $m$-stationarity, we have
\begin{equation}\label{stat}
(X_{g.U_1},\dots,X_{g.U_n})\stackrel{(d)}{=}(X^{f,m}_{at_1},\dots,X^{f,m}_{at_n})
\end{equation}
and by self-similarity,
\begin{equation}\label{self}
(X_{g.U_1},\dots,X_{g.U_n})\stackrel{(d)}{=}(\mu(g)^H X_{U_1},\dots,\mu(g)^H X_{U_n}).
\end{equation}
The result follows from (\ref{stat}) and (\ref{self}).
\fin

In the previous proof, the stationarity allows to garantee for any flow $f$ and $U\in\mathcal{A}$ the existence of $g\in G$ such that $g.U$ belongs to $f$, up to equality with respect to the law of $X$.
In that context, the $m$-stationarity definition allowing deformation of objects in $\mathcal{A}$ is the key of its special importance. 

\begin{remark}
The particular case of set-indexed fractional Brownian motion, which satisfies both properties (1) and (2) of proposition \ref{propself}, shows that projections of sifBm on any elementary flow is a self-similar one-parameter process.
Of course, this fact is already known as $m$-standard projections of sifBm are fBm. 
\end{remark}

\section{Characterization of the sifBm by stationarity and self-similarity}\label{sectcharact}
Real-parameter fractional Brownian motion is well known as the only Gaussian process satisfying the two properties of self-similarity and increment stationarity.
 
In \cite{ehem}, we proved that a set-indexed process 
$X=\left\{X_U;\;U\in\mathcal{A}\right\}$ satisfying the two following properties:
\begin{enumerate}[(i)]
\item\label{C0stat} $\mathcal{C}_0$-increment stationarity (property (\ref{oldstat}))
\item\label{selfit} self-similarity of index $H\in (0,1/2]$
\end{enumerate}
must verify, for all $U$ and $V$ in $\mathcal{A}$ such that $U\subset V$
\begin{equation*}
E\left[X_U X_V\right]=K.\left[m(U)^{2H}+m(V)^{2H}-m(V\setminus U)^{2H}\right]
\end{equation*}
where $K\in\mathbf{R}_+$.

Characterizing covariance between only comparable elements of $\mathcal{A}$, the two fractal properties (\ref{C0stat}) and (\ref{selfit}) only provide a pseudo-characterization of the sifBm.

Here, we prove that using the new definition \ref{defstat} of stationarity for a set-indexed process, we get a complete characterization of the sifBm. As we see in the proof, the statement which only consider distributional properties of set-indexed processes, relies on the characterization of the sifBm by its projections on flows (see theorem \ref{thcharactflows}).

\begin{theorem}\label{thcharact}
The sifBm $\sifbm^H$ on $(\mathcal{T},\mathcal{A},m)$ is the only $L^2$-monotone outer-continuous Gaussian set-indexed process, which is self-similar of index $H\in (0,1/2]$ and has $m$-stationary $\mathcal{C}_0$-increments.
\end{theorem}

\proof
From \cite{ehem} and proposition \ref{propsifbmstat}, we know that the sifBm is Gaussian, self-similar and has $m$-stationary $\mathcal{C}_0$-increments. 

Conversely, consider a Gaussian set-indexed process $X$, which is self-similar and has $m$-stationary $\mathcal{C}_0$-increments.
For any elementary flow $f:[a,b]\subset\mathbf{R}_+\rightarrow\mathcal{A}$, propositions \ref{statproj} and \ref{propself} imply that the standard projection of $X$ on $f$ satisfies 
\begin{itemize}
\item $X^{f,m}$ is Gaussian,
\item $X^{f,m}$ is self-similar of index $H$,
\item $X^{f,m}$ has stationary increments.
\end{itemize}
Therefore, by the well-known characterization of one-parameter fBm, $X^{f,m}$ is a fractional Brownian motion, and then
\begin{equation}\label{flow1}
\forall t\in [a,b];\quad
E\left[X_{f\circ\theta^{-1}(t)}\right]^2=t^{2H}
\end{equation}
where $\theta:t\mapsto m[f(t)]$.

Then, theorem \ref{thcharactflows} states the existence of a Borel measure $\nu$ on $\mathcal{T}$ such that $X$ is a sifBm on $(\mathcal{T},\mathcal{A},\nu)$. Consequently, the centered Gaussian process $X$ is defined by
\begin{equation*}
\forall U,V\in\mathcal{A};\quad
E\left[X_U - X_V\right]^2 = \nu(U\bigtriangleup V)^{2H},
\end{equation*}
and particularly
\begin{equation*}
\forall U\in\mathcal{A};\quad
E\left[X_U\right]^2 = \nu(U)^{2H}.
\end{equation*}
Then, according to proposition \ref{prop1}, for any elementary flow $f:[a,b]\rightarrow\mathcal{A}$, the process $\left\{X_{f\circ\psi^{-1}(t)};\;t\in [a,b]\right\}$ with $\psi:t\mapsto\nu[f(t)]$ is a one-parameter fractional Brownian motion.
This leads to
\begin{equation}\label{flow2}
\forall t\in [a,b];\quad
E\left[X_{f\circ\psi^{-1}(t)}\right]^2=t^{2H}.
\end{equation}
From (\ref{flow1}) and (\ref{flow2}), we get for any elementary flow $f:[a,b]\rightarrow\mathcal{A}$
\begin{equation*}
\forall t\in [a,b];\quad
E\left[X_{f(t)}\right]^2=m[f(t)]^{2H}=\nu[f(t)]^{2H}.
\end{equation*}
Considering a flow passing through any given $U\in\mathcal{A}$, this implies
\begin{equation*}
\forall U\in\mathcal{A};\quad
m(U)=\nu(U)
\end{equation*}
and consequently, the set-indexed process $X$ is a sifBm on $(\mathcal{T},\mathcal{A},m)$.

\fin

\begin{example}
Let us come back to example \ref{excercle}. With the same notation, the sifBm on $(\mathbb{S}_1,\mathcal{A},m)$ is the only mean-zero Gaussian process $X$ indexed by $\mathbb{S}_1$ such that the two following conditions are satisfied:
\begin{enumerate}[(i)]
\item the law of the increment process $\Delta X$ is invariant against translations along the circle;
\item $X$ is self-similar of parameter $H$, i.e.
\begin{align*}
\forall a>0;\quad
\left\{X_{a.M}; M\in\mathbb{S}_1\right\} \stackrel{(d)}{=}
\left\{a^H.X_{M}; M\in\mathbb{S}_1\right\},
\end{align*}
where $a.M$ denotes the point $M'$ of $\mathbb{S}_1$ defined by $m(\overset{\curvearrowright}{0M'}) = a.m(\overset{\curvearrowright}{0M})$.
\end{enumerate}
\end{example}

\begin{remark}
In the view of theorem \ref{thcharact}, it is natural to wonder about existence of Gaussian set-indexed processes which are self-similar of index $H\in (1/2,1)$ and have $m$-stationary $\mathcal{C}_0$-increments.
According to theorem \ref{thcharactflows}, this question is related to the existence of set-indexed processes whose standard projections on any flow are fBm of parameter $H\in (1/2,1)$.
As we saw in section \ref{sectHhalf}, the answer depends on the structure of the indexing collection $\mathcal{A}$. 
\end{remark}

\section*{Acknowledgment}
Erick Herbin would like to thank Ely Merzbach for his two kind invitations at Bar Ilan University, where most of this work was done.

\bibliographystyle{plain}
\bibliography{style}

\end{document}